\documentclass[10pt]{article}
\usepackage{amsmath}
\usepackage{amsthm}
\usepackage{amsfonts}
\usepackage{anysize}
\marginsize{2in}{2in}{1.5in}{2.375in}
\title{A note on heterogeneous decompositions into spanning trees}
\theoremstyle{plain}
\newtheorem{theorem}{Theorem}
\newtheorem{definition}{Definition}
\newtheorem{lemma}{Lemma}
\newtheorem{corollary}{Corollary}
\author{Adrian Riskin\\
Department of Mathematics\\
Mary Baldwin College\\
Staunton, VA  24401\\
ariskin@mbc.edu}
\date{}
\begin{document}
\maketitle
\begin{abstract}
In answer to a question of Eggleton, we prove that the complete multigraph on 5 vertices with
edge multiplicity 6, namely $K_{5}^{(6)}$,  has a decomposition
into 5  copies of the family of trees of order 5 and that $K_{7}^{(22)}$ has a decomposition into 
7 copies of the family of trees of order 7.  We prove something similar for $K_{2n+1}$ for 
$n \le 13$.
\end{abstract}

\section{Introduction and definitions}

Let $\mathfrak{T}(n)$ denote the family of trees of order $n$.  Let $\tau(n)=|\mathfrak{T}(n)|$.
Note that this is Sloane's A000055; see [S].
Eggleton has recently shown [E] that the complete multigraph $K_{6}^{(2)}$ of order 6 with edge multiplicity
2 has a decomposition into the elements of $\mathfrak{T}(6)$.  He asks whether $K_{5}^{(6)}$ has a 
decomposition into 5 copies of $\mathfrak{T}(5)$ and whether $K_{7}^{(22)}$ has a decomposition into
7 copies of $\mathfrak{T}(7)$.  We answer these questions in the affirmative and prove some similar results not
mentioned by Eggleton.
\begin{definition}
Let $n$, $s$, and $t$ be natural numbers.  Then
\begin{equation*}
dc_{n}(s,t)=\left\{\begin{matrix}
|s-t| \text{ if } |s-t| \le \frac{n}{2} \\
{}\\
n-|s-t| \text{ if } |s-t| \ge \frac{n}{2}
\end{matrix}\right.
\end{equation*}
\end{definition}

\begin{definition}
Let $T$ be a tree of order $2n+1$.  If it is possible to label the vertices of $T$ uniquely with the numbers 
$1, \dots, 2n+1$ so that the induced edge labels $dc_{2n+1}(\ell(u_{e}), \ell(v_{e}))$ comprise the multiset
$\{1,1,2,2,\dots,n,n\}$ then we say that $T$ is \emph{semigraceful}.
\end{definition}

\noindent This definition is of course inspired by Rosa's celebrated:

\begin{definition}
A tree of order $p$ is \emph{graceful} when the vertices can be labeled $0, \dots, p-1$ in such a way that 
the induced edge labels $|\ell (u_{e})-\ell (v_{e})|$ are distinct.  
\end{definition}

\noindent See [R] and [G] for details.  Rosa's purpose in making the definition was related to decompositions
of the complete graph into trees, although it has since acquired a vivid and independent life.  
Note that semigraceful labelings have
a certain family resemblance to the various flavors of equitable labelings, for an introduction to which see [G].

\section{Results}

The following is really more of an observation than a lemma:
\begin{lemma}
If a tree of odd order is graceful then it is semigraceful.
\end{lemma}
\noindent However there are semigraceful labelings of odd trees which are not also graceful labelings.  For instance
label the vertices of $P_{5}$, the path of order 5, with $2, 3, 1, 4, 5$.

\begin{theorem}
If $T$ is a semigraceful tree of order $2n+1$ then $K_{2n+1}^{(2)}$ is (cyclically) decomposable into $2n+1$
copies of $T$.
\end{theorem}
\begin{proof}
Draw $K_{2n+1}^{(2)}$ with its vertices evenly spaced around a circle.  Label the vertices $1, \dots, 2n+1$ in
 cyclic order.  Embed $T$ into $K_{2n+1}^{(2)}$ as directed by the semigraceful vertex labels.  This embedding
 of $T$ uses 2 edges of each possible cyclic distance in $K_{2n+1}^{(2)}$, so that when $T$ is rotated cyclically by one
 step $n$ times, each edge of $K_{2n+1}^{(2)}$ is used exactly once.
\end{proof}

\begin{corollary}
If every element of $\mathfrak{T}(2n+1)$ is semigraceful then $K_{2n+1}^{(2\tau(2n+1))}$ has a decomposition
into $2n+1$ copies of $\mathfrak{T}(2n+1)$.
\end{corollary} 

\noindent Aldred and McKay [A] have shown that every tree of order $\le 27$ is graceful, and hence 
semigraceful.  From this follows the affirmative answer to Eggleton's question about $K_{5}^{(6)}$ and
$K_{7}^{(22)}$, as well as the fact that $K_{2n+1}^{(2\tau(2n+1))}$ has such a decomposition for all 
$n \le 13$.  Now, it seems that for the most part $gcd(2n+1,\tau(2n+1))=1$ (speaking nontechnically, that is)
, and hence that $2\tau(2n+1)$ is the least possible edge multiplicity which will allow for 
such a decomposition.  However, $gcd(21, \tau(21))=gcd(21, 2144505)=3$ and $gcd(25, \tau(25))
=gcd(25, 104636890)=5$.  In these cases it is possible, at least as far as edge-counting goes, that 
$K_{21}^{(1429670)}$ is decomposable into 7 copies of $\mathfrak{T}(21)$ and that 
$K_{25}^{(41854756)}$ is decomposable into 5 copies of $\mathfrak{T}(25)$.  We conjecture that such
decompositions exist.

\section*{References}
\begin{description}

\item[[A]] Aldred, R.E.L. and McKay, B.  Graceful and harmonious labellings of trees.  
Bull. Inst. Combin. Appl.  23(1998) 69-72.

\item[[E]] Eggleton, R.B. Special heterogeneous decompositions into spanning trees.  Bull. Inst. Combin. Appl. 
 43(2005) 33-36.

\item[[G]] Gallian, J. Dynamic survey of graph labeling.  Electron. J. Combin., DS6.  
http://www.combinatorics.org/Surveys/index.html

\item[[R]] Rosa, A.  On certain valuations of the vertices of a graph.  Theory of Graphs (Int'l Symposium,
Rome, July 1966) Gordon and Breach, N.Y. 1967.

\item[[S]] Sloane, N.J.A. The on-line encyclopedia of integer sequences.  \\
http://www.research.att.com/$\sim$njas/sequences/A000055

\end{description}

\end{document}